\documentclass[reqno]{amsart}

\usepackage{amsmath}
\usepackage{amsfonts}
\usepackage{amssymb}
\usepackage{amscd}
\usepackage{amsthm}

\def\A{{\mathcal A}}

\def\1{\mathfrak{1}}
\def\0{\mathfrak{0}}

\newtheorem{lemma}{Lemma}[section]

\newtheorem{theorem}[lemma]{Theorem}
\newtheorem{proposition}[lemma]{Proposition}
\newtheorem{definition}[lemma]{Definition}
\newtheorem{remark}[lemma]{Remark}

\numberwithin{equation}{section}

\title[$C^{*}$ Estimates]{$C^{*}$ Estimates for Averaging Sums of Elements in the Thompson Group $F$}

\author{Ionut Chifan and Gabriel Picioroaga}

\address{Department of Mathematics, University of Iowa, Iowa City,USA}

\email{ichifan@math.uiowa.edu}

\address{Department of Mathematics and Computer Science, University of Southern Denmark, Odense, Denmark}

\email{gpicioro@imada.sdu.dk}



\begin{document}

\begin{abstract}
In this paper we study the non-amenability question of the Thompson Group $F$ from the 
 $C^{*}$ algebra side. Using a characterization of amenability in this framework we set about evaluating the reduced norm of the averages $\frac{1}{n}\sum x_0^ix_1x_0^{-i}$, where $x_0$ and $x_1$ are the generators of $F$ in its finite presentation. We prove that when $n$ is sufficiently large the above norm concentrates on a specific subset of $F$, easy to describe using the new normal form for elements in $F$, found by Guba and Sapir. We view this subset as the only obstruction against non-amenability. 

\end{abstract}

\maketitle

{\small \tableofcontents}

\section*{Introduction}
The Thompson group $F$ can be regarded as the group of piecewise-linear,\\ orientation-preserving 
homeomorphisms of the unit interval which have breakpoints only at dyadic points and on intervals of 
differentiability the slopes are powers of two. The group was discovered in the '60s by Richard Thompson 
and in connection with the now celebrated groups $T$ and $V$ it led to the first example of a finitely 
presented infinite simple group. Since then, these groups have received considerable applications in such 
fields as homotopy theory or operator algebras. In 1979 Geoghegan conjectured that $F$ is not amenable. 
\par In the first section we prepare some basics on the Thompson group and $C^{*}$ algebras associated with groups. We also remind a characterization of amenability in this setting. In the second section we prove the main result of the paper: for $x_0$ and $x_1$ the generators of $F$ (in the finite presentation), the averaged sum of the operators $x_0^ix_1x_0^{-i}$ converges to zero if restricted to some subspaces of $l^2(F)$. In the last section we compute the normal form of elements in $F$ for which the averaged sum may not converge to zero.

\section{Background}

\begin{definition} The Thompson group $F$ is the set of piecewise 
linear homeomorphisms from the closed unit interval $[0,1]$ to itself that are differentiable except at 
finitely many dyadic rationals and such that on intervals of differentiability the derivatives are 
powers of $2$. 
\end{definition}
\par
For a nice introduction on $F$ and its properties we refer the reader to \cite{Can}. We just remind here 
the finite and infinite presentations of $F$.
$$F=\left< x_0, x_1\mbox{ }|\mbox{ }[x_0x_1^{-1},x_0^{-1}x_1x_0]=1\mbox{, }[x_0x_1^{-1},x_0^{-2}x_1x_0^2]=1{ }\right>$$
$$F=\left< x_0, x_1, ...x_i,...|\mbox{ }x_jx_i=x_ix_{j+1}\mbox{, }i<j\mbox{ }\right>$$
It is also known that the elements of $F$ have a unique writing, the normal form. In \cite{GS}, Guba and Sapir 
found another unique way to write an element in $F$. This is what from now on in our paper will be called 
normal form. We present their result here:\\
\par $Theorem.$ (see \cite{GS}) Any element $w\in F$ can be uniquely written as a reduced word, $w=\mbox{word}(x_0,x_1)$  
such that it does not contain the following forbidden subwords:\\
1) $x_1x_0^ix_1$;\\
2) $x_1^{-1}x_0^ix_1$;\\
3) $x_1x_0^{i+1}x_1^{-1}$;\\
4) $x_1^{-1}x_0^{i+1}x_1^{-1}$,\\
for all integers $i>0$.\\
\\
Notice that $x_1^{-,+}x_0^1x_1^{-1}$ is not forbidden. 
We will apply procedures to bring certain elements of $F$ to their normal forms. For this we will need 
the following formulae (see \cite{GS}) to replace the forbidden occurences (i.e. the normal forms of the words 1,2,3,4)\\
\begin{equation}\label{e1}     
x_1x_0^ix_1=x_0^ix_1x_0^{-i-1}x_1x_0^{i+1}\quad\mbox{ }(\mbox{e.g. from  }x_{i+1}x_1=x_1x_{i+2})
\end{equation}
\begin{equation}\label{e2}     
x_1^{-1}x_0^ix_1=x_0^ix_1x_0^{-i-1}x_1^{-1}x_0^{i+1}
\end{equation}

\begin{equation}\label{e3}     
x_1x_0^{i+1}x_1^{-1}=x_0^{i+1}x_1^{-1}x_0^{-i}x_1x_0^i
\end{equation}

\begin{equation}\label{e4}     
x_1^{-1}x_0^{i+1}x_1^{-1}=x_0^{i+1}x_1^{-1}x_0^{-i}x_1^{-1}x_0^i
\end{equation}
for all $i>0$. \\
\\
This unique normal form will be combined with a characterization of amenability in the $C^*$ 
algebras realm. Next we prepare background definitions and results on group $C^{*}$ algebras. 
We also refer the reader to the books \cite{Ped} and \cite{Dav}.\\
\par Let $G$ be a countable discrete group. By $l^1(G)$ we denote the algebra of absolutely summable 
functions on $G$. The group algebra $\mathbb{C}G$, consisting of all finite sums $\sum_{g}\alpha_g\delta_g$ 
forms a dense subalgebra of $l^1(G)$. Now, the group $C^*$ algebra of G is the closure of the universal 
representation of $l^1(G)$. We denote this algebra by $C^*(G)$, its norm being defined by 
$$||f||_{C^*(G)}=\mbox{sup}\{||\pi(f)||\mbox{ }|\mbox{ }\pi\mbox{ is a *-representation of } l^1(G)\}$$
Actually we will not need this definition here, but a particular formula of the norm  when $f\in\mathbb{C}G$ (see 
\cite{Dav} for more details). Let $\mathcal{P}(G)$ denote the set of all positive definite functions 
$\phi$ on $G$ such that $\phi(e)=1$. Then, for $f=\sum_{i=1}^n\alpha_i\delta_{g_i}$ we have
$$||f||_{C^*(G)}=\mbox{sup}_{\phi\in\mathcal{P}(G)}\left(\sum_{i=1}^n\sum_{j=1}^n\alpha_i\overline{\alpha_j}\phi(g_j^{-1}g_i)\right)^{1/2}$$ 
The other $C^*$ algebra associated with the group $G$ is the so-called reduced $C^*$ algebra denoted by 
$C_r^*(G)$.\\ 
The left regular representation of $G$ on $l^2(G)$ gives rise to $C_r^{*}(G)$, 
as follows:\\
Let $l^2(G)=\left\{ \psi:G\rightarrow\mathbb{C}\mbox{ }|\mbox{}\sum_{g\in G}|\psi(g)|^2<\infty\right\}$ 
endowed with the scalar product 
$$\left<\phi,\psi\right>:=\sum_{g\in G}\phi(g)\overline{\psi(g)}$$
Notice that the Hilbert space $l^2(G)$ is generated by the countable colection of vectors 
$\left\{\delta_g\mbox{ }|\mbox{}g\in G\right\}$, where 

$$\delta_g(h)=\left\{
\begin{array}{lr}
1,&\mbox{ if }g=h \\
0,&\mbox{ otherwise} 
\end{array}\right.$$

 Also, an element $g\in G$ defines a unitary operator $\lambda(g)$, 
on $l^2(G)$  as follows: 
$\lambda(g)(\psi)(h)=\psi(g^{-1}h)$, for any $\psi\in l^2(G)$ and any $h\in G$.  Now, $C_r^{*}(G)$, 
the reduced $C^{*}$ algebra generated by $G$ is obtained by taking the norm closure in $B(l^2(G))$ of the linear span of the set 
$\left\{\lambda(g) \mbox{ }|\mbox{}g\in G\right\}$. Recall that for $A\in B(l^2(G))$ its norm is given by 
$||A||=\mbox{sup}\{||Av||\mbox{ }|\mbox{ }||v||\leq 1\}$. \\ 
The next result can be found in \cite{Ped}.\\
\par $Theorem.$ If $G$ is a locally compact group then the left regular representation $\lambda$ 
of $C^*(G)$ onto $C_r^{*}(G)$ is an isomorphism if and only if $G$ is amenable.\\
\\
We need only half of this result, which can be found in \cite{Dav} \\
\par $Theorem.$ If $G$ is a discrete amenable group then 
$$||f||_{C^*(G)}=||\lambda(f)||\mbox{    for all  }f\in\mathbb{C}G$$
\\
\\A consequence of this results is the following:\\
\par $Corollary.$  Suppose $G$ is a countable discrete group such that there exist elements 
$g_1$, $g_2$,..$g_n$ in $G$ with the property 
$$\frac{||\lambda(g_1)+\lambda(g_2)+...+\lambda(g_n)||}{n}<1$$
Then $G$ is not amenable.
\par $Proof.$  $\phi=1$ is positive definite on $G$, hence for $f=\sum_{i=1}^n\delta_{g_i}$ we obtain
$||f||_{C^*(G)}\geq n$.\\
\\
>From now on all $C^*$ norms will be reduced ones (i.e. operator norms in $B(l^2(G))$). For all preparations and proofs that 
follow we make the following convention: to not burden the notation we will write just $g$ instead 
of the operator $\lambda(g)$ and instead of the vector $\delta_g$. Any peril of confusion will be 
elliminated from the context. For example, instead of $\lambda(g)(\delta_h)=\delta_{gh}$ we will 
simply write $g(h)=gh$.\\
We are now going to make some elementary remarks about elements in $l^2(G)$. Let 
$f=\sum_{k=1}^n\alpha_kw_k$ in $l^2(G)$ such that  $w_k\in G$, for all $k$ and  
$w_k=w_l$ iff $k=l$. Then the Hilbert squared norm of $f$ is $\sum_{k=1}^n|\alpha_k|^2$. This is so  
because 
$$<w_k,w_l>=\left\{
\begin{array}{lr}
1,&\mbox{ if }k=l \\
0,&\mbox{ otherwise} 
\end{array}\right.$$

Now, if  $f=\sum_{k=1}^n\alpha_kw_k$ with possibly repeating $w$'s then we can arrange to have $f$ 
written as a finite sum $\sum_{l=1}^m\mu_lw_l$ with mutually distinct $w$'s and therefore its squared 
norm will be $\sum_{l=1}^m|\mu_l|^2$. \\
Suppose $(G_i)_{i=1}^p$ is a partition of the group $G$. For each subset $G_i$ let $H_i$ the 
Hilbert subspace it generates, i.e. the Hilbert norm closure of the linear span of $G_i$. Clearly 
$l^2(G)$ is the direct sum of the $H_i$'s. For each $i$ let $p_i$ the orthogonal projection onto $H_i$. 
We have $p_ip_j=0$ for $i\neq j$ and for $w\in G$, $p_i(w)\neq 0$ iff $w\in G_i$; in such case 
$p_i(w)=w$ (a fairly easy argument shows we cannot have both $w\in G_j$ and $w\in H_i$ for $i\neq j$). 

\section{Main result}

We will make extensive use of the following easy to check remarks. 
>From now on the shortcut "nf" stands for "normal form".\\
Let $g\in F$ written in nf and $q\in\mathbb{Z}$. Then:\\
\\
$\bullet$  $x_0^qg$  cannot contain forbidden subwords of type 1,2,3,4. Moreover, this remains 
true after reducing $x_0^qg$ (e.g. when $g$ begins with $x_0^{-,+}$).\\
\\
$\bullet$ Suppose that in its nf, $g$ begins with $x_1^{-,+}$. Then neiher $x_1^qg$ nor its reduction  
contain forbidden subwords of type 1,2,3,4.\\
\\
$\bullet$ Suppose  that in its nf, $g$ begins with $x_0^{-}$. Then $x_1^qg$ is reduced and does not 
contain forbidden subwords of type 1,2,3,4.\\
\\
\\
Next, we splitt $F$ into five disjoint subsets. Define\\
$F_1:=\{w\in F\mbox{ }|\mbox{ }w=x_0^k\mbox{, }k\geq 0\}$ and $H_1$ the closed linear span of $F_1$.\\
$F_2:=\{w\in F\mbox{ }|\mbox{ nf of }w\mbox{ begins with } x_0^kx_1^l\mbox{, }k>0\mbox{, }l\neq 0\}$ and 
$H_2$ its closed linear span.\\
$F_3:=\{w\in F\mbox{ }|\mbox{ nf of }w\mbox{ begins with } x_0^{-k}x_1^l\mbox{, }k>0\mbox{, }l\in\mathbb{Z}\}$ 
and $H_3$ its corresponding subspace. \\
$F_4:=\{w\in F\mbox{ }|\mbox{ nf of }w\mbox{ begins with } x_1^k\mbox{, }k>0\}$ and $H_4$ its corresponding 
subspace.\\
$F_5:=\{w\in F\mbox{ }|\mbox{ nf of }w\mbox{ begins with } x_1^{-k}\mbox{, }k>0\}$ and $H_5$ its 
corresponding subspace. 
The family $(F_i)_{i=1}^5$ is a partition of $F$ and therefore $l^2(F)=\bigoplus_{i=1}^5H_i$. \\
\par If there exists a constant $K<1$ and a suitable large integer $n$ such that 
$$\frac{1}{n^2}||\sum_{i=1}^nx_0^ix_1x_0^{-i}||^2<K$$
then $F$ would follow non amenable. We will prove that for large $n$ the left-hand side above 
can be replaced by  $\frac{1}{n^2}||\sum_{i=1}^nx_0^ix_1p_2x_0^{-i}||^2$, where $p_2$ is the orthogonal 
projection onto $H_2$. Thus, finding a suitable $K$ reduces to estimating the norm (in $B(l^2(F))$ this time) 
of the averaged sum at vectors $v\in H_2$.  

\begin{proposition}\label{p2} Let $p_i$ the orthogonal projection onto $H_i$ and $p:=p_1+p_3+p_4+p_5$. 
Then we have the following estimate:
$$||\sum_{i=1}^nx_0^ix_1x_0^{-i}||^2\leq 8||\sum_{i=1}^nx_0^ix_1p_1x_0^{-i}||^2+
8||\sum_{i=1}^nx_0^ix_1p_3x_0^{-i}||^2+$$ 
$$+8||\sum_{i=1}^nx_0^ix_1p_4x_0^{-i}||^2+8||\sum_{i=1}^nx_0^ix_1p_5x_0^{-i}||^2+$$ 
$$+2||\sum_{i=1}^nx_0^ix_1p_2x_0^{-i}||\mbox{ }||\sum_{j=1}^nx_0^jx_1px_0^{-j}||+||\sum_{i=1}^nx_0^ix_1p_2x_0^{-i}||^2$$
\end{proposition}
\begin{proof} The reason we left the squared norm containing $p_2$ at the end is that we did not wanted 
it be multiplied by too large a constant (larger than 1, actually). Also, the theorem below will 
be more illuminating. Let us proceed with the proof. Notice first 
$x_1=x_1p_1+x_1p_3+x_1p_4+x_1p_5+x_1p_2$, because of the partition of $F$.  
Also, for operators $A$ and $B$ in some $B(H)$, recall the following inequality:
$$||A+B||^2=||(A+B)^*(A+B)||\leq ||A||^2 + ||B||^2 + 2||A||\mbox{ }||B||$$
Apply this inequality for $A=\sum_{i=1}^nx_0^ix_1px_0^{-i}$ and $B=\sum_{i=1}^nx_0^ix_1p_2x_0^{-i}$. 
For majorizing $||A||^2$  use three times the inequality $||C+D||^2\leq 2||C||^2 + 2||D||^2$. 
\end{proof}
\begin{theorem} Let $p_i$ as above and $p=p_1+p_3+p_4+p_5$. We have:  
$$\mbox{a) }\lim_{n\rightarrow\infty}\frac{1}{n^2}||\sum_{i=1}^nx_0^ix_1p_1x_0^{-i}||^2=0$$     
$$\mbox{b) }\lim_{n\rightarrow\infty}\frac{1}{n^2}||\sum_{i=1}^nx_0^ix_1p_3x_0^{-i}||^2=0$$
$$\mbox{c) }\lim_{n\rightarrow\infty}\frac{1}{n^2}||\sum_{i=1}^nx_0^ix_1p_4x_0^{-i}||^2=0$$
$$\mbox{d) }\lim_{n\rightarrow\infty}\frac{1}{n^2}||\sum_{i=1}^nx_0^ix_1p_5x_0^{-i}||^2=0$$
$$\mbox{e) }\lim_{n\rightarrow\infty}\frac{1}{n^2}||\sum_{i=1}^nx_0^ix_1p_2x_0^{-i}||\mbox{ }||\sum_{j=1}^nx_0^jx_1px_0^{-j}||=0$$
\end{theorem}
\begin{proof} a) Let $w\in l^2(F)$, $w=\sum_{k=1}^m\alpha_kw_k$, $w_k\in F$ such that 
$||w||^2=\sum |\alpha|^2\leq 1$
We have 
$$||\sum_{i=1}^nx_0^ix_1p_1x_0^{-i}(w)||^2=\sum_{i,j}\sum_{k,l}\alpha_k\overline{\alpha_l}<x_0^ix_1p_1(x_0^{-i}w_k), x_0^jx_1p_1(x_0^{-j}w_l)>$$
In the above sum the scalar products for which $p_1(x_0^{-i}w_k)=0$ or $p_1(x_0^{-j}w_l)=0$ do not count. 
In general, for $g\in F$, $p_1(x_0^{-i}g)\neq 0$ iff $x_0^{-i}g\in F_1$ that is $g=x_0^h$ with $h\geq i$. 
In the sum above, a non zero scalar product would have to be equal to 1 and  would have to have the form:
$<x_0^ix_1x_0^{h_k}, x_0^jx_1x_0^{h_l}>$, where $h_k\geq i$, $h_l\geq j$. This happens iff
$x_0^ix_1x_0^{h_k}=x_0^jx_1x_0^{h_l}$. By the uniqueness of the nf we must have $i=j$ and $w_k=w_l$ i.e. 
the norm above can be majorized: 
$$||\sum_{i=1}^nx_0^ix_1p_1x_0^{-i}(w)||^2\leq \sum_{i=1}^n\sum_{k=1}^m|\alpha_k|^2\leq n$$
In general, for $A\in B(H)$, $||A||=\mbox{sup}\{||A(w)||\mbox{ : }||w||\leq 1\}$, but in our situation it is not hard  
 to see that we can take the sup over vectors $w$ as above. In conclusion a) follows. We will do the same 
 for the other limits, the idea being to use the uniqueness of the nf.\\
\\
b) We have
$$||\sum_{i=1}^nx_0^ix_1p_3x_0^{-i}(w)||^2=\sum_{i,j}\sum_{k,l}\alpha_k\overline{\alpha_l}<x_0^ix_1p_3(x_0^{-i}w_k), x_0^jx_1p_3(x_0^{-j}w_l)>$$
As above, for a non zero scalar product 
$$<x_0^ix_1x_0^{-i}w_k, x_0^jx_1x_0^{-j}w_l>=1\mbox{     and     }x_0^{-i}w_k\in F_3\mbox{ , }x_0^{-j}w_l\in F_3$$
Let us discuss the restriction $x_0^{-i}g\in F_3$ where $g\in F$ is written in its nf. We want to prove 
that $x_1x_0^{-i}g$ cannot contain forbidden subwords of type 1), 2), 3) or 4) and after reduction the 
nf of $x_0^ix_1x_0^{-i}g$ begins with $x_0^ix_1$. We will use the remarks at the beginning of the section.
\par -- if $g\in F_1$ then $g=x_0^k$. Because $x_0^{-i}g\in F_3$ we must have $k<i$. This implies 
that $x_0^ix_1x_0^{-i}g=x_0^ix_1x_0^{-i+k}$ and the left-hand side is the nf we wanted.\\
\par -- if $g\in F_2$ then $g=x_0^kx_1^l...$ in nf, $k>0$, $l\neq 0$. We must have $k<i$ (otherwise, using again 
the remarks above $x_0^{-i}g=x_0^{-i+k}x_1^l...\notin F_3$). We conclude that the nf of 
$x_0^ix_1x_0^{-i}g$ is $x_0^ix_1x_0^{-i+k}x_1^l...$, as we wanted.\\
\par -- if  $g\in F_3$ then $x_1x_0^{-i}g$ cannot contain forbidden subwords, and after reducing this 
word the letter $x_1$ will still maintain its first position. Therefore, left multypling by $x_0^i$ 
will produce a word that begins with $x^ix_1$ in nf, which is what we wanted.\\
\par -- if $g\in F_4\cup F_5$. In this case $x_0^{-i}g$ is already in nf and so does 
$x_0^ix_1x_0^{-i}g$.\\

\par In conclusion, for a non zero scalar product to appear, it is required that \\
$x_0^ix_1x_0^{-i}w_k= x_0^jx_1x_0^{-j}w_l$. By the uniqueness of nf and cases above we infer $i=j$. 
Simplifying the last equality we obtain also $w_k=w_l$, situation which we agreed to happen iff $k=l$. 
Now we can conclude b). \\
\\
c) As above, the discussion comes down to $p_4(x_0^{-i}g)\neq 0$. Hence $g=x_0^ix_1^h...$, $h>0$. 
We easily infer (using the remarks again) that $x_0^ix_1p_4(x_0^{-i}g)=x_0^ix_1^{1+h}...$, where 
the last element is in nf. Therefore $x_0^ix_1x_0^{-i}w_k= x_0^jx_1x_0^{-j}w_l$ implies 
$x_0^ix_1^{1+h_k}...=x_0^jx_1^{1+h_l}...$. The last equality of nf's implies again $i=j$ and $k=l$.\\
\\
d) As in c) $g=x_0^ix_1^h...$, but this time $h<0$. The discussion is similar to the one in c) 
if $h<-1$: we obtain that $x_0^ix_1x_0^{-i}g$ begins with $x_0^i$ in nf. If $h=-1$, the $x_1$ letter 
in the middle will cancel, probably affecting the first position $x_0^i$. In this situation is possible 
to have $i\neq j$ for some of the scalar products, eventhough the nf's settle immediately. 
We will use the following trick (which works for $p_4$ as well):
\par Assume $b\in B(l^2(F))$ such that $x_0^ibx_0^{-i}b^*=0$ for all 
$i=\overline{1\mbox{,...}n}$. Then $$||\sum_{i=1}^nx_0^ibx_0^{-i}||^2\leq n||b||^2$$
Notice that $x_0$ is unitary. Also, if $AB^*=0$ then $||A+B||^2\leq ||A||^2 + ||B||^2$. We have 
$$||\sum_{i=1}^nx_0^ibx_0^{-i}||^2=||x_0[b + x_0^{-1}(\sum_{i=2}^nx_0^ibx_0^{-i})x_0]x_0^{-1}||^2
\leq ||b||^2 + ||\sum_{i=1}^{n-1}x_0^ibx_0^{-i}||^2$$
Inductively, we obtain the desired estimate. Now, we apply it for $b=x_1p_5$. Of course, we have 
to make sure  $x_0^ix_1p_5x_0^{-i}p_5x_1^{-1}=0$. Actually $p_5x_0^{-i}p_5x_1^{-1}(w)=0$ for any $w\in F$ such that $p_5x_0^{-i}p_5(w)=0$. This last equality is easy to prove: first, $w\in F_5$, otherwise $p_5(w)=0$. Left multiplying a $F_5$-word by $x_0^{-}$ does not produce forbidden subwords, therefore $x_0^{-i}p_5(w)=x_0^{-i}w\in F_3$. Projecting in $F_5$ gives $p_5x_0^{-i}p_5(w)=0$.
\par In conclusion  $$||\sum_{i=1}^nx_0^ix_1p_5x_0^{-i}||^2\leq n||x_1||^2\mbox{ }||p_5||^2\leq n$$ 
and d) follows.\\
\\
e) Follows from a), b), c) and d).
\end{proof} 
\begin{remark}Let us notice here that the same treatment for 
$||\sum_{i=1}^nx_0^ix_1p_2x_0^{-i}(w)||$ does not work.  For $w=\sum\alpha_kw_k\in H_2$ it is possible to have $x_0^ix_1x_0^{-i}w_k=x_0^jx_1x_0^{-j}w_l$ for $i\neq j$. Eventhough we can further "minimize" such occurences (see next section, Proposition \ref{min}), it is  possible for a word in $F_2$ that begins with $x_0^p$, $p$ large, to land in $F_2$, under the action of $x_0^ix_1x_0^{-i}$. As a consequence we should not expect the limit of $p_2$ averages be zero. Actually, an upper bound strictly less than $1$ would be just enough, however we have not been able to do this. It would seem necessary to efficiently count the pairs $(i,j)$ and $(k,l)$ for which $x_0^ix_1x_0^{-i}w_k=x_0^jx_1x_0^{-j}w_l$. Thus finding  the normal forms of the elements $x_0^ix_1x_0^{-i}w$ for $w\in F_2$ might be helpful. This is done in the next section. 
\end{remark}

\section{Normal forms in the subset $F_2$}
The proof of the next result consists of a straightforward computation with the aid of 
formulae (\ref{e1}),..,(\ref{e4}). 
\begin{proposition}\label{p1} Let $H>0$ and $L\neq 0$ integers such that the elements $x_1x_0^Hx_1^L$ and 
$x_1^{-1}x_0^Hx_1^L$ of $F$, contain forbidden subwords. Then their normal forms are:\\
\begin{equation}\label{f1}
x_1^{-1}x_0^Hx_1^L=x_0^Hx_1^Lx_0^{-H-L}x_1^{-1}x_0^{H+L}\mbox{, if }H+L>0
\end{equation}
\begin{equation}\label{f2}
x_1^{-1}x_0^Hx_1^L=x_0^Hx_1^{-H+1}x_0^{-1}x_1^{-1}x_0x_1^{H+L-1}\mbox{, if }H+L\leq 0
\end{equation}
\begin{equation}\label{f3}
x_1x_0^Hx_1^L=x_0^Hx_1^Lx_0^{-H-L}x_1x_0^{H+L}\mbox{, if }H+L>0
\end{equation}
\begin{equation}\label{f4}
x_1x_0^Hx_1^L=x_0^Hx_1^{-H+1}x_0^{-1}x_1x_0x_1^{H+L-1}\mbox{, if }H+L\leq 0
\end{equation}
\end{proposition}

\par Notice that (\ref{f2}) and (\ref{f4}) do not make sense for $H=1$ (e.g. $x_1^{-,+}x_0x_1^L$ is not forbidden if $1+L\leq 0$).
\par Before describing the normal forms of the elements 
$x_0^ix_1x_0^{-i}w$ for $w\in F_2$ we prove the result mentioned at the end of the previous section.
\begin{proposition}\label{min}
Let $w_1$ and $w_2$ two distinct elements of the Thompson's group $F$. Then there may exist only one pair $(i,j)$ of positive integers such that 
$$x_0^ix_1x_0^{-i}w_1=x_0^jx_1x_0^{-j}w_2$$
\end{proposition}
\begin{proof}Notice first that for any such pair, $i\neq j$ as $w_1\neq w_2$. Assume there are two pairs $i\neq j$ and $k\neq l$ that satisfy the equation above. Solving for $w_1$ and $w_2$ we obtain
$$x_0^ix_1^{-1}x_0^{j-i}x_1x_0^{-j}=
x_0^kx_1^{-1}x_0^{l-k}x_1x_0^{-l}$$
If $j-i<0$ and $l-k<0$ then in both sides of the equality we have normal forms. By uniqueness, we get $i=k$ and $j=l$.\\
If $j-i<0$ and $l-k>0$ then the left-hand side is a nf. We work out the right-hand side with the aid of (\ref{e2}) and obtain :
$$x_0^ix_1^{-1}x_0^{j-i}x_1x_0^{-j}=
x_0^lx_1x_0^{-l+k-1}x_1^{-1}x_0^{-k+1}$$
This is impossible as both terms must be normal forms and the second occurence in both expressions ($x_1\neq x_1^{-1}$) does not match. \\
If $i-j>0$ and $l-k<0$ then a similar argument implies a contradiction.\\
If $i-j>0$ and $l-k>0$ then after applying (\ref{e2}) we obtain the normal forms:
$$x_0^jx_1x_0^{-j+i-1}x_1^{-1}x_0^{-i+1}=x_0^lx_1x_0^{-l+k-1}x_1^{-1}x_0^{-k+1}$$
Again, by uniqueness we get $j=l$ and $i=k$. 
\end{proof}
We begin now to describe the normal form of the elements
$$iw:=x_0^ix_1x_0^{-i}w$$
A typical element in $F_2$ (in nf) looks like 
$$w=x_0^hx_1^lx_0^{p_1}x_1^{q_1}\mbox{...}x_0^{p_m}x_1^{q_m}x_0^{p_{m+1}}$$
subject to restrictions (in order to achieve nf of $F_2$):\\
$h>0$, $l\neq 0$, $m\geq 0$, $p_{m+1}\in\mathbb{Z}$, $0\neq p_i\leq 1$ (with $p_i=1$ only if $q_i<0$),\\ 
$0\neq q_i$ for all $i\leq m$. 
\par In order to bring $iw$ to nf we will apply the formulae in Proposition \ref{p1}. We will be interested only in the case $h>i$, otherwise there are no forbidden subwords. Remark also that if $h-i=1$ and $l<0$ then the nf is again obtained immediately. For the other situations we will obtain six  possible types of normal forms. We prefer to not summarize it in a proposition, but rather display every nf obtained along the computations.
\par The first forbidden occurence in $iw$ is $x_1x_0^{h-i}x_1^l$. We may apply either (\ref{f3}) or (\ref{f4}). After that, in order to remove forbidden subwords we may only apply (\ref{f3}) or (\ref{f4}).\\
\\ I. Let us see what we get if we apply (\ref{f3}) $t$ times, $t<m+1$:
\begin{equation}\label{nf1}
iw=x_0^hx_1^lx_0^{p_1}x_1^{q_1}...x_0^{p_{t-1}}x_1^{q_{t-1}}x_0^{-c}x_1x_0^{c+p_t}x_1^{q_t}...x_0^{p_{m+1}}
\end{equation}
where $c:=h-i+l+\sum_{k=1}^{t-1}(p_k+q_k)$ subject to the following restrictions (that appear in order to apply (\ref{f3}) up to step $t$):
\begin{equation}\label{nf10}
h-i+l+\sum_{k=1}^{\nu}(p_k+q_k)>0, \mbox{ for all }0\leq\nu<t
\end{equation}
\begin{equation}\label{nf11}
h-i+l+\sum_{k=1}^{\nu-1}(p_k+q_k)+p_{\nu}\geq 1, \mbox{ for all }0\leq\nu<t\mbox{ with "=" only if } q_{\nu}>0 
\end{equation}
At step $t$ the only possible forbidden occurence that makes the procedure go forward (to the right in the sequence) is $x_1x_0^{c+p_t}x_1^{q_t}$. If we hit a value $t$ such that 
\begin{equation}\label{nf12}
 c+p_t=h-i+l+\sum_{k=1}^{t-1}(p_k+q_k)+p_t\leq 1\mbox{ with "=" only if }q_t<0
\end{equation}
then we cannot apply (\ref{f3}) or (\ref{f4}) anymore.  However, $iw$ may not be reduced yet.\\
\\ I.1 If $c+p_t\neq 0$ and (\ref{nf12}) holds then $iw$ is reduced and (\ref{nf1}) is a normal form.\\
\\ I.2  If (\ref{nf10}), (\ref{nf11}) hold but $c+p_t=0$ then (\ref{nf1}) can be written
\begin{equation}\label{nf2}
iw=x_0^hx_1^lx_0^{p_1}x_1^{q_1}...x_0^{p_{t-1}}x_1^{q_{t-1}}x_0^{p_t}x_1^{1+q_t}x_0^{p_{t+1}}x_1^{q_{t+1}}...x_0^{p_{m+1}}
\end{equation}
Let us prove that (\ref{nf2}) is a nf provided $1+q_t\neq 0$. The subword $x_1^{1+q_t}x_0^{p_{t+1}}x_1^{q_{t+1}}$ is not forbidden because $w$ is in nf and restrictions apply to $p_{t+1}$ and $q_{t+1}$. Also $p_t=-c<0$ thus $x_1^{q_{t-1}}x_0^{p_t}x_1^{1+q_t}$ is not forbidden. Hence (\ref{nf2}) is a nf.\\
\\ I.3 If $1+q_t=0$ in (\ref{nf2}) then $iw$ can be written
\begin{equation}\label{nf3}
iw=x_0^hx_1^lx_0^{p_1}x_1^{q_1}...x_0^{p_{t-1}}x_1^{q_{t-1}}x_0^{p_t+p_{t+1}}x_1^{q_{t+1}}...x_0^{p_{m+1}}
\end{equation}
We argue that if $p_t+p_{t+1}\neq 0$ then (\ref{nf3}) is a nf. The only possible forbidden occurence is $x_1^{q_{t-1}}x_0^{p_t+p_{t+1}}x_1^{q_{t+1}}$. Remark that if $p_t+p_{t+1}<0$ then we do not have a forbidden occurence. Assume by contradiction $p_t+p_{t+1}>0$. Because $p_t=-c<0$ we get $p_{t+1}>0$. We also get $p_t\leq -1$ as all powers are integers. The restrictions for $w$ imply $p_{t+1}=1$. We obtain the  contradiction\\
$0<p_t+p_{t+1}=1+p_t\leq 0$.\\
\\ I.4 If $p_t+p_{t+1}=0$ in (\ref{nf3}) we will show that the procedure of bringing $iw$ to nf will stop. Recall all restrictions so far: $c>0$, $c+p_t=0$ and $1+q_t=0$. With $p_t+p_{t+1}=0$, $iw$ is now written
\begin{equation}\label{nf4}
iw=x_0^hx_1^lx_0^{p_1}x_1^{q_1}...x_0^{p_{t-1}}x_1^{q_{t-1}+q_{t+1}}x_0^{p_{t+2}}x_1^{q_{t+2}}...x_0^{p_{m+1}}
\end{equation}
We claim that the right-hand side of (\ref{nf4}) is a nf, i.e. suffices to show that $q_{t-1}+q_{t+1}<0$ (recall that $x_1x_0x_1^{-}$ is not forbidden). As above, we have $p_t<0$ and from $p_t+p_{t+1}=0$ we necessarilly get $p_t=-1$ and $p_{t+1}=1$. The restrictions on the nf of $w$ show $q_{t+1}<0$. We will also prove $q_{t-1}<0$ and then we're done. By contradiction  suppose $q_{t-1}>0$ (in a nf only non-zero powers show up). Put $\nu=t-1$ in the restriction (\ref{nf11}):
$$h-i+l+\sum_{k=1}^{t-2}(p_k+q_k)+p_{t-1}\geq 1$$
Adding up $q_{t-1}>0$, the inequality becomes:  
$$h-i+l+\sum_{k=1}^{t-1}(p_k+q_k)>1$$
Taking into account $p_t=-1$ we rewrite:
$$h-i+l+\sum_{k=1}^{t-1}(p_k+q_k)+p_t>0$$
However, this last inequality contradicts $c+p_t=0$. In conclusion, if we apply (\ref{f3}) repeatedly  then we obtain four types of normal forms of the element $iw$.\\
 \\ II. If (\ref{f4}) is to be applied from the start ( to get rid of the first forbidden occurence $x_1x_0^{h-i}x_1^l$ in $iw$) then $h-i+l\leq 0$ and nf of $iw$ settles immediately. This is actually a special case of III below:
\begin{equation}
iw=x_0^hx_1^{-h+i+1}x_0^{-1}x_1x_0x_1^{h-i+l-1}x_0^{p_1}x_1^{q_1}...x_0^{p_{m+1}}
\end{equation}
Again, $x_1x_0x_1^{h-i+l-1}$ is not forbidden as $h-i+l-1<0$\\
\\III. Apply (\ref{f3}) a couple of times such that at step $t-1$ conditions to apply (\ref{f4}) are fulfilled. We prove that nf of $iw$ settles down in at most two steps after applying (\ref{f4}).By looking at (\ref{nf1}) and its restrictions we obtain
\begin{equation}\label{nf5} 
iw=x_0^hx_1^lx_0^{p_1}x_1^{q_1}...x_0^{p_{t-2}}x_1^{q_{t-2}}x_0^{p_{t-1}}x_1^{-d}x_0^{-1}x_1x_0x_1^{d+q_{t-1}}x_0^{p_t}x_1^{q_t}...x_0^{p_{m+1}}
\end{equation}
where $d:=h-i+l+\sum_{k=1}^{t-2}(p_k+q_k)+p_{t-1}-1$ and the following restrictions apply:
\begin{equation}
h-i+l+\sum_{k=1}^{\nu-1}(p_k+q_k)+p_{\nu}\geq 1, \mbox{ for all }0\leq\nu<t\mbox{ with "=" only if } q_{\nu}>0 
\end{equation}
\begin{equation}\label{nf20}
h-i+l+\sum_{k=1}^{\nu-1}(p_k+q_k)> 0, \mbox{ for all  }0\leq\nu<t 
\end{equation}
\begin{equation}\label{nf21}
h-i+l+\sum_{k=1}^{t-1}(p_k+q_k)\leq 0 
\end{equation}
The last inequality comes from the condition $H+L\leq 0$ ($H$, $L$ correspond to the forbidden subword in  (\ref{nf1})), needed to apply (\ref{f4}). \\
\\
 III.1 If $d$ is non-zero we claim that the right-hand side of (\ref{nf5}) is a nf. Notice  $d+q_{t-1}<0$ by (\ref{nf21}). Hence, forbidden occurences could only appear at $x_1^{q_{t-2}}x_0^{p_{t-1}}x_1^{-d}$. This happens only if $p_{t-1}=1$ and $-d>0$, which would contradict (\ref{nf20}) for $\nu=t-1$. In conclusion, (\ref{nf5}) is a normal form. \\
 \\
 III.2 If $d=0$ in (\ref{nf5}) we prove that nf of $iw$ is 
\begin{equation}\label{nf6}
iw=x_0^hx_1^lx_0^{p_1}x_1^{q_1}...x_0^{p_{t-2}}x_1^{q_{t-2}}x_0^{p_{t-1}-1}x_1x_0x_1^{q_{t-1}}x_0^{p_t}x_1^{q_t}...x_0^{p_{m+1}}
\end{equation}
It suffices to prove that the procedure cannot go further left, i.e. $p_{t-1}-1<0$, when $d=0$. If both numbers are zero then inequality (\ref{nf20}) would be violated for $\nu=t-1$. Also, there are no forbidden occurences: the only possible spot for such a subword is $x_1x_0x_1^{q_{t-1}}$. However, we must have $d+q_{t-1}<0$. Hence, (\ref{nf6}) is the last type of normal form an element $iw$ can have.

\end{document}